\documentclass{amsart}
\usepackage{amssymb, amsthm, amsmath}

\newtheorem{theorem}{Theorem}[section]
\newtheorem{lemma}[theorem]{Lemma}

\theoremstyle{definition}

\theoremstyle{remark}
\newtheorem{remark}[theorem]{Remark}
\newtheorem{example}[theorem]{Example}

\def\hilb{\mathop{\mathrm{Hilb}}}
\def\Hom{\mathop{\mathrm{Hom}}}
\def\ann{\mathop{\mathrm{ann}}}
\def\HK{\mathrm{HK}}
\def\m{\mathfrak{m}}
\def\a{\mathfrak{a}}
\def\b{\mathfrak{b}}
\def\N{\mathbb{N}}
\def\Z{\mathbb{Z}}
\def\Q{\mathbb{Q}}

\begin{document}
\title{The F-signature of an affine semigroup ring}
\author{Anurag K. Singh}
\address{
School of Mathematics \\
686 Cherry St. \\
Georgia Institute of Technology \\
Atlanta, GA 30332-0160, USA.
E-mail: {\tt singh@math.gatech.edu}}

\thanks{This material is based upon work supported in part by the National
Science Foundation under grants DMS-0070268 and DMS-0300600.}

\dedicatory{Dedicated to Professor Kei-ichi Watanabe on the occasion of his
sixtieth birthday}

\begin{abstract}
We prove that the F-signature of an affine semigroup ring of positive
characteristic is always a rational number, and describe a method for computing
this number. We use this method to determine the F-signature of Segre products
of polynomial rings, and of Veronese subrings of polynomial rings. Our
technique involves expressing the F-signature of an affine semigroup ring as
the difference of the Hilbert-Kunz multiplicities of two monomial ideals, and
then using Watanabe's result that these Hilbert-Kunz multiplicities are
rational numbers.
\end{abstract}

\maketitle

\section{Introduction}

Let $(R,\m)$ be a Cohen-Macaulay local or graded ring of characteristic $p>0$,
such that the residue field $R/\m$ is perfect. We assume that $R$ is reduced
and F-finite. Throughout $q$ shall denote a power of $p$, i.e., $q=p^e$ for
$e \in \mathbb N$. Let
$$
R^{1/q} \approx R^{a_q} \oplus M_q
$$
where $M_q$ is an $R$-module with no free summands. The number $a_q$ is
unchanged when we replace $R$ by its $\m$-adic completion, and hence is
well-defined by the Krull-Schmidt theorem. In \cite{HL} Huneke and Leuschke
define the {\it F-signature}\/ of $R$ as
$$
s(R) = \lim_{q \to \infty} \frac{a_q}{q^{\dim R}},
$$
provided this limit exists. In this note we study the F-signature of normal
monomial rings, and our main result is:

\begin{theorem} \label{main}
Let $K$ be a perfect field of positive characteristic, and $R$ be a normal
subring of a polynomial ring $K[x_1, \dots, x_n]$ which is generated, as a
$K$-algebra, by monomials in the variables $x_1, \dots, x_n$. Then the
F-signature $s(R)$ exists and is a positive rational number.

Moreover, $s(R)$ depends only on the semigroup of monomials generating $R$ and
not on the characteristic of the perfect field $K$. 
\end{theorem}

We also develop a general method for computing $s(R)$ for monomial rings, and 
use it to determine the F-signature of Segre products of polynomial rings, and
of Veronese subrings of polynomial rings.

In general, it seems reasonable to conjecture that the limit $s(R)$ exists and
is a rational number. Huneke and Leuschke proved that the limit exists if $R$
is a Gorenstein ring, \cite[Theorem 11]{HL}. They also proved that a ring $R$
is weakly F-regular whenever the limit is positive, and this was extended by
Aberbach and Leuschke in \cite{AL}:

\begin{theorem}[Huneke-Leuschke, Aberbach-Leuschke]
Let $(R,\m)$ be an F-finite reduced Cohen-Macaulay ring of characteristic
$p>0$. Then $R$ is strongly F-regular if and only if
$$
\limsup_{q\to\infty} \frac{a_q}{q^{\dim R}}>0.
$$
\end{theorem}

Further results on the existence of the F-signature are obtained by Aberbach
and Enescu in the recent preprint \cite{AE}. Also, the work of Watanabe and
Yoshida \cite{WY} and Yao \cite{Yao} is closely related to the questions
studied here.

We mention that a graded $R$-module decomposition of $R^{1/q}$ was used by
Peskine-Szpiro, Hartshorne, and Hochster, to construct small Cohen-Macaulay
modules for $R$ in the case where $R$ is an $\N$-graded ring of dimension three,
finitely generated over a field $R_0$ of characteristic $p>0$, see
\cite[\S 5 F]{melwitt}. The relationship between the $R$-module decomposition
of $R^{1/q}$ and the singularities of $R$ was investigated in the paper
\cite{SV}.

\section{Semigroup rings}

The semigroup of nonnegative integers will be denoted by $\N$. Let $x_1, \dots,
x_n$ be variables over a field $K$. By a {\it monomial}\/ in the variables
$x_1, \dots, x_n$, we will mean an element $x_1^{h_1} \cdots x_n^{h_n} \in
K[x_1, \dots, x_n]$ where $h_i \in \N$. We frequently switch between semigroups
of monomials in $x_1, \dots, x_n$ and subsemigroups of ${\N}^n$, where we
identify a monomial $x_1^{h_1} \cdots x_n^{h_n}$ with $(h_1,\dots,h_n) \in
{\N}^n$. A semigroup $M$ of monomials is {\it normal}\/ if it is finitely
generated, and whenever $a,b$ and $c$ are monomials in $M$ such that $ab^k=c^k$
for some positive integer $k$, then there exists a monomial $\alpha \in M$ with
$\alpha^k=a$. It is well-known that a semigroup $M$ of monomials is normal if
and only if the subring $K[M] \subseteq K[x_1, \dots, x_n]$ is a normal ring,
see \cite[Proposition 1]{meltoric}. 

A semigroup $M$ of monomials is {\it full}\/ if whenever $a,b$ and $c$ are 
monomials such that $ab=c$ and $b,c \in M$, then $a \in M$. By
\cite[Proposition 1]{meltoric}, a normal semigroup of monomials is isomorphic
(as a semigroup) to a full semigroup of monomials in a possibly different set
of variables.

\begin{lemma} \label{sequence}
Let $A=K[x_1, \dots, x_n]$ be a polynomial ring over a field $K$, and
$R \subseteq A$ be a subring generated by a full semigroup of monomials. Let
$\m$ denote the homogeneous maximal ideal of $R$, and assume that $R$ contains
a monomial $\mu$ in which each variable $x_i$ occurs with positive exponent.
For positive integers $t$, let $\a_t$ denote the ideal of $R$ generated by
monomials in $R$ which do not divide $\mu^t$.

\begin{enumerate}
\item The ideals $\a_t$ are irreducible and $\m$-primary, and the image of
$\mu^t$ generates the socle of the ring $R/\a_t$.

\item The ideals $\a_t$ form a non-increasing sequence $\a_1 \supseteq \a_2
\supseteq \a_3 \supseteq \dots$ which is cofinal with the sequence $\m
\supseteq \m^2 \supseteq \m^3 \supseteq \dots$.

\item Let $M$ be a finitely generated $R$-module with no free summands. Then
$\mu^t M \subseteq \a_t M$ for all $t \gg 0$.

\item Let $K$ be a perfect field of characteristic $p>0$, and 
$R^{1/q} \approx R^{a_q} \oplus M_q$ be an $R$-module decomposition
of $R^{1/q}$ where $M_q$ has no free summands. Then 
$$
a_q = \ell \left( \frac{R}{\a_t^{[q]} :_R \mu^{tq}} \right)
\quad \text{for all} \quad t \gg 0.
$$
\end{enumerate}
\end{lemma}

\begin{proof}

(1) It suffices to consider $t=1$ and $\a=\a_1$. Every non-constant monomial in
$R$ has a suitably high power which does not divide $\mu$, so $\a$ is
$\m$-primary. If $\alpha \in R$ is any monomial of positive degree, then
$\alpha \mu \in \a$, and so $\m \subseteq \a :_R \mu$. Also $\mu \notin \a$, so
we conclude that $\a :_R \mu = \m$. Since $\a$ is a monomial ideal, the socle
of $R/\a$ is spanned by the images of some monomials. If $\theta \in R$ is a
monomial whose image is a nonzero element of the socle of $R/\a$, then $\mu =
\beta\theta$ for a monomial $\beta \in R$. If $\beta \in \m$ then $\mu \in \m
\theta \subseteq \a$, a contradiction. Consequently we must have $\beta=1$,
i.e., $\theta=\mu$.

(2) Since each $x_i$ occurs in $\mu \in R$ with positive exponent and $R$ is
generated by a full semigroup of monomials, we see that 
$$
\a_t \subseteq (x_1^{t+1}, \dots, x_n^{t+1})A \cap R.
$$
It follows that $\{ \a_t \}_{t \in {\N}}$ is cofinal with
$\{ \m^t \}_{t \in {\N}}$.

(3) For an arbitrary element $m \in M$, consider the homomorphism $\phi: R \to
M$ given by $r \mapsto rm$. Since the module $M$ has no free summands, $\phi$
is not a split homomorphism. By \cite[Remark 2]{melnagoya}, there exists $t_0
\in \N$ such that $\mu^{t_0} m \in \a_{t_0} M$, equivalently, such that the
induced map
$$
\overline{\phi}_{t_0}: R/\a_{t_0} \to M/\a_{t_0} M
$$
is not injective. If $\overline{\phi}_t: R/\a_t \to M/\a_t M$ is injective for 
some $t \ge t_0$, then it splits since $R/\a_t$ is a Gorenstein ring of 
dimension zero; however this implies that 
$$
\overline{\phi}_{t_0} : 
 R/\a_t \otimes_{R/\a_t} R/\a_{t_0} \to M/\a_t M \otimes_{R/\a_t} R/\a_{t_0}
$$
splits as well, a contradiction. Consequently
$\overline{\phi_t}(\overline{\mu^t}) = 0$, i.e., $\mu^t m \in \a_t M$ for all
$t \gg t_0$. The module $M$ is finitely generated, and so
$\mu^t M \subseteq \a_t M$ for all $t \gg 0$. 

(4) For any ideal $\b \subseteq R$, we have
$$
\frac{R^{1/q}}{\b R^{1/q}} \cong \left( \frac{R}{\b R} \right)^{a_q} \oplus 
\frac{M_q}{\b M_q},
$$
and so
$$
\ell \left( \frac{R}{\b^{[q]}} \right) 
 = \ell \left( \frac{R^{1/q}}{\b R^{1/q}} \right) 
 = a_q \ell \left( \frac{R}{\b} \right) 
 + \ell \left( \frac{M_q}{\b M_q} \right). 
$$
Using this for the ideals $\a_t$ and $\a_t + \mu^tR$ and taking the difference, 
we get
\begin{multline*}
a_q \left[ \ell \left( \frac{R}{\a_t} \right) 
 - \ell \left( \frac{R}{\a_t + \mu^t R} \right) \right]
 + \ell \left( \frac{M_q}{\a_t M_q} \right)
 - \ell \left( \frac{M_q}{\a_t M_q + \mu^t M_q} \right) \\
= \ell \left( \frac{R}{\a_t^{[q]}} \right) 
- \ell \left( \frac{R}{\a_t^{[q]} + \mu^{tq} R} \right) 
= \ell \left( \frac{R}{\a_t^{[q]} :_R \mu^{tq}} \right) 
\end{multline*}
By (3) $\mu^t M_q \subseteq \a_t M_q$ for all $t \gg 0$, and the result follows.
\end{proof}

\begin{lemma} \label{independent}
Let $K$ be a perfect field of characteristic $p>0$, and $R$ be a subring of
$A=K[x_1,\dots, x_n]$ generated by a full semigroup of monomials with the
property that for every $i$ with $1 \le i \le n$, there exists a monomial $a_i
\in A$ in the variables $x_1, \dots, \widehat{x_i}, \dots, x_n$ such that
$a_i/x_i = \mu_i/\eta_i$ for monomials $\mu_i, \eta_i \in R$. Let $\mu_0 \in R$
be a monomial in which each $x_i$ occurs with positive exponent, and set $\mu =
\mu_0\mu_1 \cdots \mu_n$. For $t \ge 1$, let $\a_t$ be the ideal of $R$
generated by monomials in $R$ which do not divide $\mu^t$. Then, for every
prime power $q=p^e$ and integer $t \ge 1$, we have
$$
\a_t^{[q]} :_R \mu^{tq} = \m_A^{[q]} \cap R 
$$
where $\m_A=(x_1, \dots, x_n)A$ is the maximal ideal of $A$. If 
$R^{1/q} \approx R^{a_q} \oplus M_q$ is an $R$-module decomposition of
$R^{1/q}$ where $M_q$ has no free summands, then
$$
a_q = \ell \left( \frac{R}{\a_t^{[q]} :_R \mu^{tq}} \right)
= \ell \left( \frac{R}{\m_A^{[q]} \cap R} \right)
\quad \text{for all} \quad q=p^e \quad \text{and} \quad t \ge 1.
$$
\end{lemma}

\begin{proof} 
By Lemma \ref{sequence}\,$(4)$, it suffices to prove that 
$$
\a_t^{[q]} :_R \mu^{tq} = \m_A^{[q]} \cap R \quad \text{for all} \quad q=p^e
\quad \text{and} \quad t \ge 1.
$$
Given a monomial $r \in \a_t^{[q]} :_R \mu^{tq}$, there exists a monomial $\eta
\in R$ which does not divide $\mu^t$ for which $r \mu^{tq} \in \eta^q R$. Since
$\mu^t/\eta$ is an element of the fraction field of $R$ which is not in $R$, we
must have $\mu^t/\eta \notin A$ and so $\eta A :_A \mu^t \subseteq \m_A$.
Taking Frobenius powers over the regular ring $A$, we get
$$
\eta^q A :_A \mu^{tq} \subseteq \m_A^{[q]}
$$
and hence $r \in \m_A^{[q]} \cap R$. This shows that
$\a_t^{[q]} :_R \mu^{tq} \subseteq \m_A^{[q]} \cap R$.

For the reverse inclusion, consider a monomial $bx_i^q \in R$ where $b \in A$.
Then
$$
bx_i^q \mu^{tq} = b a_i^q \left( \frac{\mu^t \eta_i}{\mu_i} \right)^q 
$$
where $ba_i^q$ and $\mu^t \eta_i / \mu_i$ are elements of $R$. It remains to
verify that $\mu^t \eta_i/\mu_i \in \a_t$, i.e., that it does not divide $\mu^t$
in $R$. Since
$$
\frac{\mu^t}{\mu^t \eta_i/\mu_i} = \frac{a_i}{x_i},
$$
this follows immediately.
\end{proof}

\begin{lemma} \label{embed}
Let $R'$ be a normal monomial subring of a polynomial ring over a field $K$. 
Then $R'$ is isomorphic to a subring $R$ of a polynomial ring
$A=K[x_1,\dots,x_n]$ where $R$ is generated by a full semigroup of
monomials, and for every $1 \le i \le n$, there exists a monomial $a_i \in A$
in the variables $x_1, \dots, \widehat{x_i}, \dots, x_n$, for which $a_i/x_i$
is an element of the fraction field of $R$.
\end{lemma}

\begin{proof}
Let $M \subseteq {\N}^r$ be the subsemigroup corresponding to 
$R' \subseteq K[y_1,\dots,y_r]$. Let $W \subseteq {\Q}^r$ denote the $\Q$-vector
space spanned by $M$, and $W^*=\Hom_{\Q}(W,{\Q})$ be its dual vector space.
Then
$$
U= \{ w^* \in W^* : w^*(m) \ge 0 \quad \text{for all} \quad m \in M \}
$$
is a finite intersection of half-spaces in $W^*$. Let $w_1^*,\dots,w_n^* \in U$
be a minimal ${\Q}_+$-generating set for $U$, where ${\Q}_+$ denotes the
nonnegative rationals. Replacing each $w_i^*$ by a suitable positive multiple,
we may ensure that $w_i^*(m) \in {\N}$ for all $m \in M$, and also that
$w_i^*(M) \nsubseteq a \Z$ for any integer $a \ge 2$. It is established in
\cite[\S 2]{meltoric} that the map $T: W \to {\Q}^n$ given by
$$
T=(w_1^*, \dots, w_n^*) 
$$
takes $M$ to an isomorphic copy $T(M) \subseteq {\N}^n$, which is a full
subsemigroup of ${\N}^n$. Let $R \subseteq A=K[x_1,\dots,x_n]$ be the monomial
subring corresponding to $T(M) \subseteq {\N}^n$.

Fix $i$ with $1 \le i \le n$. Since $w_i^*(M) \nsubseteq a \Z$ for any integer
$a \ge 2$, the fraction field of $R$ contains an element
$x_1^{h_1}\cdots x_n^{h_n}$ such that $h_1,\dots,h_n \in \Z$ and $h_i = -1$.
Also, there exists $m \in M$ such that $w_i^*(m) = 0$ and $w_j^*(m) \neq 0$ for
all $j \neq i$. Consequently $R$ contains a monomial
$\alpha = x_1^{s_1} \cdots x_n^{s_n}$ with $s_i = 0$ and $s_j > 0$ for all
$j \neq i$. For a suitably large integer $t \ge 1$, the element
$$
x_1^{h_1} \cdots x_n^{h_n} \alpha^t = a_i/x_i
$$
belongs to the fraction field of $R$ where $a_i\in A$ is a monomial in the
variables $x_1,\dots,\widehat{x_i},\dots,x_n$.
\end{proof}

\begin{proof}[Proof of Theorem \ref{main}]
By Lemma \ref{embed}, we may assume that $R$ is a monomial subring of 
$A=K[x_1,\dots,x_n]$ satisfying the hypotheses of Lemma \ref{independent}. For
the choice of $\mu$ as in Lemma \ref{independent}, the ideals
$\a_t^{[q]} :_R \mu^{tq}$ do not depend on $t \in \mathbb N$. Setting $\a = \a_1$
we get
$$
a_q = \ell \left( \frac{R}{\a^{[q]} :_R \mu^q} \right) 
= \ell \left( \frac{R}{\a^{[q]}} \right) 
- \ell \left( \frac{R}{\a^{[q]} + \mu^q R} \right),
$$
i.e., $a_q$, as a function of $q=p^e$, is a difference of two Hilbert-Kunz
functions. Let $d = \dim R$. By \cite{monsky} the limits
$$
e_{\HK}(\a)=\lim_{q\to\infty} \frac{1}{q^d} \ell \left(\frac{R}{\a^{[q]}}\right) 
\quad \text{ and } \quad 
e_{\HK}(\a+\mu R)=\lim_{q\to\infty} \frac{1}{q^d} 
 \ell \left(\frac{R}{\a^{[q]}+\mu^q R}\right)
$$
exist, and by \cite{watanabe} they are rational numbers. Consequently the limit
$$
\lim_{q \to \infty} \frac{a_q}{q^d} = e_{\HK}(\a) - e_{\HK}(\a+\mu R)
$$
exists and is a rational number. The ring $R$ is F-regular, so the positivity 
of $s(R)$ follows from the main result of \cite{AL}; as an alternative proof,
we point out that $\mu\notin\a^*$, and consequently
$e_{\HK}(\a) > e_{\HK}(\a+\mu R)$ by \cite[Theorem 8.17]{HHjams}.

By \cite{watanabe} the Hilbert-Kunz multiplicities $e_{\HK}(\a)$ and
$e_{\HK}(\a+\mu R)$ do not depend on the characteristic of the field $K$,
and so the same is true for $s(R)$.
\end{proof}

\begin{remark}
Let $(R,\m,K)$ be a local or graded ring of characteristic $p>0$, and let $\eta
\in E_R(K)$ be a generator of the socle of the injective hull of $K$. In
\cite{WY} Watanabe and Yoshida define the {\it minimal relative Hilbert-Kunz
multiplicity} of $R$ to be
$$
m_{\HK}(R) = \liminf_{e \to \infty} \frac{\ell(R/\ann_R(F^e(\eta)))}{p^{de}},
$$
where $d = \dim R$. They compute $m_{\HK}(R)$ in the case $R$ is the Segre
product of polynomial rings, \cite[Theorem 5.8]{WY}. Their work is closely
related to our computation of $s(R)$ in the example below.
\end{remark}

\section{Examples}

\begin{example}
Let $K$ be a perfect field of positive characteristic, and consider integers
$r,s \ge 2$. Let $R$ be the Segre product of the polynomial rings
$K[x_1,\dots,x_r]$ and $K[y_1,\dots,y_s]$, i.e., $R$ is subring of
$A=K[x_1,\dots,x_r,y_1,\dots,y_s]$ generated over $K$ be the monomials $x_iy_j$
for $1 \le i \le r$ and $1 \le j \le s$. It is well-known that $R$ is
isomorphic to the determinantal ring obtained by killing the size two minors of
an $r \times s$ matrix of indeterminates, and that the dimension of the ring
$R$ is $d=r+s-1$. Lemma \ref{independent} enables us to compute not just the
F-signature $s(R)$, but also a closed-form expression for the numbers $a_q$.

The rings $R \subseteq A$ satisfy the hypotheses of Lemma \ref{independent},
and so
$$
a_q = \ell \left( \frac{R}{ \m_A^{[q]} \cap R } \right) 
 = \ell \left( \frac{K[x_1,\dots,x_r]}{(x_1^q,\dots,x_r^q)}
 \ \# \ \frac{K[y_1,\dots,y_s]}{(y_1^q,\dots,y_s^q)} \right),
$$
where $\#$ denotes the Segre product. The Hilbert-Poincar\'e series of these 
rings are
$$
\hilb \left( \frac{K[x_1,\dots,x_r]}{(x_1^q,\dots,x_r^q)},u \right) 
 = \frac{(1-u^q)^r}{(1-u)^r} \quad \text{and} \quad 
\hilb \left( \frac{K[y_1,\dots,y_s]}{(y_1^q,\dots,y_s^q)},v \right) 
 = \frac{(1-v^q)^s}{(1-v)^s}
$$
and so $a_q$ is the sum of the coefficients of $u^iv^i$ in the polynomial
$$
\frac{(1-u^q)^r}{(1-u)^r} \frac{(1-v^q)^s}{(1-v)^s} \in {\mathbb Z}[u,v].
$$
Therefore $a_q$ equals the constant term of the Laurent polynomial
$$
\frac{(1-u^q)^r}{(1-u)^r} \frac{(1-u^{-q})^s}{(1-u^{-1})^s} 
 = \frac{ u^s (1-u^q)^{r+s} }{ u^{sq} (1-u)^{r+s} } 
 \in {\mathbb Z}[u, u^{-1}],
$$
and hence the coefficient of $u^{s(q-1)}$ in 
$$
\frac{(1-u^q)^{r+s}}{(1-u)^{r+s}} = 
\left[ \sum_{i = 0}^{r+s} (-1)^i \binom{r+s}{i}u^{iq} \right] 
\left[ \sum_{n \ge 0} \binom{d+n}{d} u^n \right] .
$$ 
Consequently we get
$$
a_q = \sum_{i = 0}^s (-1)^i \binom{r+s}{i} \binom{d + s(q-1) - iq}{d}
= \sum_{i = 0}^s (-1)^i \binom{d+1}{i} \binom{q(s-i) + d - s}{d}
$$
where we follows the convention that $\binom{m}{n} = 0$ unless $0 \le n \le m$.
This shows that the F-signature of $R$ is
$$
s(R) = \lim_{q \to \infty} \frac{a_q}{q^d} = 
\frac{1}{d \, !} \sum_{i = 0}^s (-1)^i \binom{d+1}{i} (s-i)^d.
$$
We point out that $s(R) = A(d,s)/d\,!$ where the numbers
$$
A(d,s) = \sum_{i = 0}^s (-1)^i \binom{d+1}{i} (s-i)^d
$$
are the {\it Eulerian numbers}, i.e., the number of permutations of $d$ objects
with $s-1$ {\it descents}; more precisely, $A(d,s)$ is the number of 
permutations $\pi = a_1 a_2 \cdots a_d \in S_d$ whose {\it descent set}
$$
D(\pi) = \{ i : a_i > a_{i+1} \}
$$
has cardinality $s-1$, see \cite[\S 1.3]{stanley}. These numbers satisfy 
the recursion
$$
A(d,s)=sA(d-1,s)+(d-s+1) A(d-1,s-1) \quad \text{ where } \quad A(1,1)=1. 
$$
\end{example}

\begin{example}
Let $K$ be a perfect field of positive characteristic. For integers $n\ge 1$
and $d \ge 2$, let $R$ be the $n$\,th Veronese subring of the polynomial ring
$A=K[x_1,\dots,x_d]$, i.e., $R$ is subring of $A$ which is generated, as a
$K$-algebra, by the monomials of degree $n$. In the case $d=2$ and $p \nmid n$,
the F-signature of $R$ is $s(R)=1/n$, as worked out in \cite[Example 17]{HL}.

It is readily seen that the rings $R \subseteq A$ satisfy the hypotheses of 
Lemma \ref{independent}, and therefore 
$$
a_q = \ell \left( \frac{R}{ \m_A^{[q]} \cap R } \right).
$$
Consequently $a_q$ equals the sum of the coefficients of $1, t^n, t^{2n}, \dots$
in
$$
\hilb \left( \frac{K[x_1,\dots,x_d]}{(x_1^q,\dots,x_d^q)}, t \right) 
= \frac{(1-t^q)^d}{(1-t)^d} = (1+t+t^2+ \dots + t^{q-1})^d.
$$
Let $f(m)$ be the sum of the coefficients of powers of $t^n$ in
$$
(1+t+t^2+ \dots + t^{m-1})^d.
$$
A routine computation using, for example, induction on $d$, gives us
$f(n) = n^{d-1}$, and it follows that
$$
f(kn)= k^df(n) = k^d n^{d-1}.
$$
To obtain bounds for $a_q=f(q)$, choose integers $k_i$ with
$k_1n \le q \le k_2n$ where $0 \le |q-k_in| \le n-1$. Then
$f(k_1n) \le f(q) \le f(k_2n)$, and hence
$$
\left(\frac{q-n+1}{n}\right)^d n^{d-1} \le k_1^d n^{d-1} \le a_q
\le k_2^d n^{d-1} \le \left(\frac{q+n-1}{n}\right)^d n^{d-1}
$$
Consequently
$$
a_q = \frac{q^d}{n} + O(q^{d-1}),
$$
and $s(R)=1/n$ .
\end{example}

\section*{Acknowledgments}

I would like to thank Ezra Miller for several useful discussions.


\begin{thebibliography}{Ho2}

\bibitem[AE]{AE} I. M. Aberbach and F. Enescu, {\em When does the F-signature
exist?\/} preprint.

\bibitem[AL]{AL} I. M. Aberbach and G. J. Leuschke, {\em The F-signature and
strong F-regularity}, Math. Res. Lett. {\bf 10} (2003), 51--56.

\bibitem[Ho1]{meltoric} M. Hochster, {\em Rings of invariants of tori,
Cohen-Macaulay rings generated by monomials, and polytopes}, Ann. of Math. (2)
{\bf 96} (1972), 318--337.

\bibitem[Ho2]{melnagoya} M. Hochster, {\em Contracted ideals from integral
extensions of regular rings}, Nagoya Math. J. {\bf 51} (1973), 25--43.

\bibitem[Ho3]{melwitt} M. Hochster, {\em Big Cohen-Macaulay modules and
algebras and embeddability in rings of Witt vectors}, Conference on Commutative
Algebra--1975 (Queen's Univ., Kingston, Ont., 1975), pp. 106--195. Queen's
Papers on Pure and Applied Math., No. 42, Queen's Univ., Kingston, Ont., 1975.

\bibitem[HH]{HHjams} M. Hochster and C. Huneke, {\em Tight closure, invariant
theory, and the Brian\c con-Skoda theorem}, J. Amer. Math. Soc. {\bf 3} (1990),
31--116.

\bibitem[HL]{HL} C. Huneke and G. J. Leuschke, {\em Two theorems about maximal
Cohen-Macaulay modules}, Math. Ann. {\bf 324} (2002), 391--404.

\bibitem[Mo]{monsky} P. Monsky, {\em The Hilbert-Kunz function}, Math. Ann.
{\bf 263} (1983), 43--49.

\bibitem[SV]{SV} K. E. Smith and M. Van den Bergh, {\em Simplicity of rings of
differential operators in prime characteristic}, Proc. London Math. Soc. (3)
{\bf 75} (1997), 32--62.

\bibitem[St]{stanley} R. Stanley, {\em Enumerative combinatorics, Volume I},
Wadsworth \& Brooks/Cole Advanced Books \& Software, Monterey, CA, 1986. 

\bibitem[Wa]{watanabe} K.-i. Watanabe, {\em Hilbert-Kunz multiplicity of toric
rings}, Proc. Inst. Natural Sci., Nihon Univ., {\bf 35} (2000), 173--177.

\bibitem[WY]{WY} K.-i. Watanabe and K.-i. Yoshida {\em Minimal relative
Hilbert-Kunz multiplicity}, Illinois J. Math. {\bf 428} (2004), 273--294.

\bibitem[Ya]{Yao} Y. Yao, {\em Observations on the F-signature of local rings
of characteristic p}, preprint.

\end{thebibliography}
\end{document}